\documentclass[10pt,final]{article}

\usepackage{graphicx}
\usepackage{amssymb, amsthm}

\newcommand{\calS}{\mathcal{S}}
\newcommand{\E}{\mathbb{E}}

\newcommand{\C}{\mathbb{C}}
\newcommand{\R}{\mathbb{R}}
\newcommand{\norm}[1]{\left\|#1\right\|}

\newcommand{\lb}{\left(}
\newcommand{\lsb}{\left[}
\newcommand{\lcb}{\left\{}

\newcommand{\lab}{\left\langle}
\newcommand{\lAb}{\left\langle\langle}
\newcommand{\rb}{\right)}
\newcommand{\rsb}{\right]}
\newcommand{\rcb}{\right\}}

\newcommand{\rab}{\right\rangle}
\newcommand{\rAb}{\right\rangle\rangle}
\newcommand{\erf}{\mathrm{erf}}

\newtheorem{theorem}{Theorem}[section]
\newtheorem{definition}[theorem]{Definition}

\newtheorem{lemma}[theorem]{Lemma}
\newtheorem{remark}[theorem]{Remark}
\newtheorem{example}[theorem]{Example}
\newtheorem{proposition}[theorem]{Proposition}
\newenvironment{pf}{{\bf Proof.}}{{\unskip\nobreak\hfil\penalty50\hskip2em\hbox{}\nobreak\hfil\rule{2mm}{2mm}\parfillskip=0pt \finalhyphendemerits=0 \par}\vskip 2ex}

\begin{document}
\begin{center}
\begin{LARGE}Parameter Estimation from Occupation Times\end{LARGE} \\
\begin{small}
\end{small}
\end{center}
\vskip 10mm
\begin{center}
\begin{large}W.~Bock \end{large}   \\
\textit{\begin{small}bock@mathematik.uni-kl.de \\
\end{small} }
\begin{small}
\textit{Functional Analysis and Stochastic Analysis Group, \\
Department of Mathematics, \\
University of Kaiserslautern, 67653 Kaiserslautern, Germany \\}
\end{small}
\vskip 10mm
\begin{large} T.~G\"otz \end{large} \\   
\textit{\begin{small}goetz@uni-koblenz.de \\
\end{small} }
\begin{small}
\textit{Mathematical Institute,\\
University of Koblenz, Universit{\"a}tsstr.~1, 56070 Koblenz, Germany\\}
\end{small}
\vskip 10mm
\begin{large} M.~Grothaus \end{large} \\   
\textit{\begin{small}grothaus@mathematik.uni-kl.de \\
\end{small} }
\begin{small}
\textit{Functional Analysis and Stochastic Analysis Group, \\
Department of Mathematics, \\
University of Kaiserslautern, 67653 Kaiserslautern, Germany \\}
\end{small}
\vskip 10mm
\begin{large} U.~P.~Liyanage \end{large} \\   
\textit{\begin{small}liyanage@mathematik.uni-kl.de \\
\end{small} }
\begin{small}
\textit{Technomathematics Group, \\
Department of Mathematics, \\
University of Kaiserslautern, 67653 Kaiserslautern, Germany \\}
\end{small}
\end{center}
\pagebreak
\begin{abstract}
We derive an equation to compute directly the expected occupation time of the centered Ornstein--Uhlenbeck process. This allows us to identify the parameters of the Ornstein--Uhlenbeck process for available occupation times via a standard least squares minimization. To test the method, we generate occupation times via Monte--Carlo simulations and recover the parameters with the above mentioned procedure.
\end{abstract}

\section{Introduction}

Nonwoven materials or fleece are webs of long flexible fibers that are used for composite materials, e.g.~filters, as well as in the hygiene and textile industries.  They are produced in melt--spinning operations: hundreds of individual endless fibers are obtained by the continuous extrusion of a molten polymer through narrow nozzles that are densely and equidistantly placed in a row at a spinning beam. The viscous or viscoelastic fibers are stretched and spun until they solidify due to cooling air streams. Before the elastic fibers lay down on a moving conveyor belt to form a web, they become entangled and form loops due to the highly turbulent air flows. The homogeneity and load capacity of the fiber web are the most important textile properties for quality assessment of industrial nonwoven fabrics. The optimization and control of the fleece quality require modeling and simulation of fiber dynamics and lay--down. Available data to judge the quality, at least on the industrial scale, are usually the mass per unit area of the fleece. 

A stochastic model for the fiber deposition in the nonwoven production was proposed and analyzed in Ref.~\cite{BGKMW08,GKMW07}. Its core is a stochastic Ornstein--Uhlenbeck process for the random motion of the fiber. The aim of this paper is to determine the parameters of the Ornstein--Uhlenbeck process from available mass per unit area data, i.e.~the occupation time in mathematical terms. For the sake of simplicity, we focus on a one--dimensional version of the Ornstein--Uhlenbeck process. 

The paper is organized as follows: In Section~\ref{S:Model} we introduce the Ornstein--Uhlenbeck process as a prototypic model for the fiber deposition. Section~\ref{S:Theory} is devoted to the derivation of the expectation value for the occupation times. An algorithm to estimate the parameters in the Ornstein--Uhlenbeck process from available occupation times is presented along with numerical experiments in Section~\ref{S:Numerics}. Finally, we draw some conclusions and give an outlook to open questions.
 
\section{Model}
\label{S:Model}

As a prototypic model for the fiber deposition, we consider the general one--dimensional Ornstein--Uhlenbeck process
\begin{equation}
   \label{E:SDE_OU}
    d U_t = \lambda(\mu-U_t) dt + \sigma dW_t, \quad U_0 = U(0)\in \R\;.
\end{equation}
The equilibrium $\mu\in \R$ and the stiffness $\lambda>0$ govern the deterministic part whereas the diffusion parameter $\sigma>0$ and a standard Wiener process (or Brownian motion) $W_t$ contribute the stochastic part. For the sake of simplicity, we consider mainly the centered process $X_t$ satisfying
\begin{equation}
    \label{E:SDE_OU_SIM}
    d X_t = -\lambda X_t dt + \sigma dW_t, \quad X_0 = 0\;.
\end{equation}

The real-valued random variable $X_t$ models the deposition point of an individual fiber on the fleece.  If we follow the random variable over a time interval $[0,T]$ for $T>0$, we obtain the path of an individual fiber. To introduce the mathematical analog of the mass per unit area we need the following definition.

\begin{definition}[Occupation time]
Let $T>0$ and consider an interval $[a,b]\subset\R$, where $a=-\infty$ or $b=\infty$ are allowed. The \emph{occupation time} $M_{T,[a,b]}$ is defined as
$$
    M_{T,[a,b]}(X_t) := \int_0^T \mathbf{1}_{[a,b]}(X_t) dt = \int_0^T\int_a^b \delta_0(X_t-x)\, dx\, dt\;.
$$
Here, $\mathbf{1}_{[a,b]}$ denotes the indicator function of the interval $[a,b]$ and $\delta_0(X_t-x)$ is the Donsker's delta function introduced in Definition~\ref{D:Donsker}, below.
\end{definition}

\begin{remark}
The occupation time is a random variable itself. It models the time, the random process spends inside the spatial interval $[a,b]$ during the time interval $[0,T]$. In terms of our physical model for the nonwoven production, the occupation time can be interpreted as the mass of fiber material deposited inside the interval $[a,b]$, i.e.~the mass per unit area of the final fleece. This quantity is easily accessible even on the scale of industrial production and hence it will serve as the input to our parameter estimation problem.
\end{remark}

In the next chapter, we will present tools from white noise analysis to derive the expectation of the occupation time for the centered Ornstein--Uhlenbeck process $X_t$. Although it is possible to derive the results by classical stochastic analysis methods, we use a white noise approach to generalize the concepts also to higher dimensions, where one can give a rigorous meaning to multidimensional Donskers Delta functions as a white noise distribution, in later research. Moreover in future work an extension to more complicated processes (e.g. with fractional noise term) is planed. Thereafter we show, how to estimate the parameters $\lambda$, $\sigma$ of the process from available data for the occupation times.

\section{Theory}
\label{S:Theory}
We start by considering the Gel'fand triple $S(\R) \subset L^2(\R) \subset S'(\R)$, where $S(\R)$ denotes the Schwartz space of rapidly decreasing smooth functions, $L^2(\R)$ the Hilbert space of real--valued square integrable (equivalence classes of) functions on $\R$ w.r.t.~Lebesgue measure and $S'(\R)$ the topological dual of $S(\R)$, i.e.~the space of tempered distributions. This particular choice is the usual one in white noise analysis \cite{HKPS93}. By $\lab f ,\omega  \rab$ we denote the duality pairing between $\omega \in S'(\R)$ and $f\in S(\R)$, an extension of the standard inner product on $L^2(\R)$ in the sense of a Gel'fand triple.

Next, we want to introduce a probability measure on the space $S'(\R)$. Therefore, we consider the $\sigma$--algebra $\mathcal{B}(S'(\R))$ generated by the cylinder sets $\lcb \lab f, \cdot\rab:\ f\in S(\R)\rcb$. The white noise measure $\mu$ on $(S'(\R),\mathcal{B})$ is given via Minlos' theorem \cite{BK95,Hi80,HKPS93} by its characteristic function $C$
$$
   \int_{S'(\R)} \exp(i \lab  f,\omega  \rab) \, d\mu(\omega) = \exp \lb -\frac{1}{2} |f|^2 \rb = C(f)
$$
for $f \in S(\R)$. 

\begin{remark} 
The Hilbert space of complex--valued square--integrable functions w.r.t.~this measure $\mu$ is denoted by $L^2(\mu)=L^2(S'(\R),\mathcal{B},\mu)$. For $f,g\in S(\R)$ we have the isometry
$$
   \int_{S'(\R)} \lab  f,\omega  \rab\,  \lab  g,\omega  \rab \, d\mu(\omega) 
   	= \int_{\R} f(s) g(s) \, ds\;.
$$
Thus, this result can also be extended to $f,g\in L^2(\R)$ in the sense of an $L^2(\mu)$--limit. Hence, within the above formalism, a version of a standard Wiener process can be written as $W_t=\lab  \mathbf{1}_{[0,t)}, \cdot \rab$, for $t>0$ and $W_0=0$.
\end{remark}

To treat the occupation time of the Ornstein--Uhlenbeck process in the white noise framework, we need the space of Hida distributions $(\mathcal{S})'$.

The above introduced space $L^2(\mu)$ serves as the central space of the Gel'fand triple $(\mathcal{S})\subset L^2(\mu)\subset (\mathcal{S})'$, where $(\mathcal{S})$ denotes the space of Hida test functions. The dual pairing of $\Phi\in (\mathcal{S})'$ with $\varphi\in (\mathcal{S})$ is denoted by $\lAb \varphi, \Phi\rAb$. For a detailed description of the construction of the Hida triple we refer to Ref.~\cite{HKPS93}.

\begin{example}
For a function $f\in S(\R)$, the exponential $ \exp(i\lab f,\cdot\rab)$ is an element of $(\mathcal{S})$.
\end{example}

We will characterize Hida distributions with the help of the $T$--transform and $U$--functionals.

\begin{definition}[$T$--transform] $\,$\\
The $T$--transform of a Hida distribution $\Phi\in (\calS)'$ is defined as
$$
    T(\Phi)(f) := \lAb \Phi, \exp(i\lab f,\cdot\rab) \rAb\;,
$$
where $f\in S(\R)$.

Since $ 1\in(\calS)$, the expectation of a Hida distribution $\Phi\in (\calS)'$ can be defined by
$$
    \E_{\mu}(\Phi) := \lAb 1,\Phi\rAb = T(\Phi)(0)\;.
$$
\end{definition}

\begin{definition}[$U$-functional] $\,$\\
We call $F:S(\R)\to \C$ a $U$-functional, if
\begin{enumerate}
\item For all $f,g\in S(\R)$, the mapping $\R\ni x\mapsto F(g+xf)\in \C$ is analytic and hence has an entire extension to $\C$.
\item There exist constants $0\le K,C<\infty$ and a continuous norm $\norm{\cdot}$ on $S(\R)$ such that 
$$
    |F(z\xi)| \le K \exp(c |z|^2\norm{\xi}^2)\;,
$$
for all $z\in \C$ and all $\xi\in S(\R)$.
\end{enumerate}
\end{definition}

The proof of the following equivalence theorem can be found in Ref.~\cite{HKPS93}.
\begin{theorem}
\label{T:Char}
A mapping $F:S(\R)\to \C$ is the $T$--transform of a unique element in $(\calS)'$, if and only if $F$ is a $U$--functional.
\end{theorem}

\begin{example}
In the sense of a limit in $(\calS)'$ we can define the white noise process as
$$
     \omega(t) := \lab \delta_t,\omega  \rab \in (\calS)'\;,
$$
where $\delta_t$ denotes the Dirac delta in $t>0$. This process can be considered as the time derivative of the Wiener process $W_t(\omega)=\langle {\bf 1}_{[0,t)},\omega \rangle$ in the sense of Hida distributions.
\end{example}

The next result follows from Theorem \ref{T:Char} and concerns integration of a family of Hida distributions, see~Ref.~\cite{HKPS93,KLPSW96,PS91}.
\begin{theorem}
\label{T:Intcor}
Let $(\Lambda, \mathcal{A}, \nu)$ be a measure space and $\lambda \mapsto \Phi(\lambda)$ a mapping from $\Lambda$ to $(\calS)'$. We assume that the $T$--transform $T(\Phi(\lambda))$ satisfies the following conditions:
\begin{enumerate}
\item The mapping $\lambda \mapsto T(\Phi(\lambda))(f)$ is measurable for all $f \in S(\R)$.
\item There exists a continu\-ous norm $\norm{\cdot}$ on $S(\R)$ and functions \\
$C \in L^{\infty}(\mathcal{A}, \nu)$ and $D \in L^1(\mathcal{A}, \nu)$ integrable with respect to $\nu$ such that 
$$
   |T(\Phi(\lambda))(zf)| \leq D(\lambda)\cdot \exp(C(\lambda) |z|^2 \norm{f}^2)\;,   
$$
for all $f \in S(\R)$, $z\in \C$ .
\end{enumerate}
Then it holds in the sense of Bochner integration in a suitable sub--Hilbert space of $(\calS)'$, that the integral of the familiy of Hida distributions is itself a Hida distribution, i.e.~$\displaystyle \int_{\Lambda} \Phi(\lambda) \, d\nu(\lambda) \in (\calS)'$ and the $T$--transform interchanges with the integration
$$
   T\lb \int_{\Lambda} \Phi(\lambda) \, d\nu(\lambda) \rb =
   	\int_{\Lambda} T(\Phi(\lambda)) \, d\nu(\lambda)
$$
\end{theorem}

Based on the above theorem, we introduce the following Hida distribution.
\begin{definition}[Donsker]
\label{D:Donsker} 
We define Donsker's delta at $x \in \R$ corresponding to $\eta \in L^2(\R)$ by
$$
   \delta_x(\lab \eta,\cdot \rab) := 
   	\frac{1}{2\pi} \int_{\R} \exp(i \lambda (\lab \eta,\cdot \rab -x)) \, d \lambda
$$
in the sense of Bochner integration\cite{HKPS93,LLSW94,W95}. Its $T$--transform in $f \in S(\R)$ is given by
$$
   T(\delta_x(\lab  \eta,\cdot \rab)(f) 
   	= \frac{1}{\sqrt{2\pi \lab \eta, \eta\rab}} \exp\left( -\frac{1}{2\lab \eta,\eta \rab}(i\lab \eta,f \rab - x)^2 -\frac{1}{2}\lab f,f\rab \right)\;.
$$
\end{definition}
\bigskip
Coming back to the Ornstein--Uhlenbeck process, we note that we can write it in the framework of white noise analysis as 
\begin{equation}
    U_t(\omega) = (U_0-\mu) \exp(-\lambda t) +\mu + \lab \sigma \exp\lsb\lambda(\cdot-t)\rsb\mathbf{1}_{[0,t)}(\cdot), \omega \rab \in L^2(\mu)\;,
\end{equation}
for $t\ge 0$. This can be seen as follows:\\
Clearly $U_t$ is a Gaussian random variable with expectation
\begin{eqnarray}
   \mathbb{E}_{\mu}(U_t) &= \int_{S'(\R)} U_t  \, d\mu \nonumber
		= U_0 \exp(-\lambda t) +\mu(1-\exp(\lambda t)) \nonumber
\end{eqnarray}		
and covariance
\begin{eqnarray}
   {\rm Cov}(U_t,U_{\tau}) &=& \mathbb{E}_{\mu} \big( \lb U_t- \mathbb{E}_{\mu}(U_t)\rb \cdot \lb U_{\tau}-\mathbb{E}_{\mu}(U_{\tau})\rb\big) \nonumber\\
	&=& \int_{S'(\R)} \lab  \sigma \exp \lsb \lambda(\cdot-t)\rsb \mathbf{1}_{[0,t)}(\cdot),\omega \rab \nonumber \\
	&&\quad\quad \quad \times \lab \sigma \exp\lsb \lambda(\cdot-\tau)\rsb \mathbf{1}_{[0,\tau)}(\cdot),\omega \rab \, d\mu \nonumber\\
	&=& \frac{\sigma^2}{2\lambda} \big( \exp(-\lambda\cdot |t-\tau|)-\exp(- \lambda (t+\tau))\big) \nonumber.
\end{eqnarray}
Thus, by uniqueness $U_t$ is the Ornstein--Uhlenbeck process solving the corresponding SDE~(\ref{E:SDE_OU}).

In the special case $U_0=\mu=0$ of the centered Ornstein--Uhlenbeck process $X_t=\lab \eta_t , \cdot \rab$, where $s \mapsto \eta_t(s) = \sigma \exp \lsb \lambda(s -t)\rsb \mathbf{1}_{[0,t)}(s)\, \in L^2(\R)$, we obtain that the $T$--transform of the corresponding Donsker's delta at $x\in \R$ is given by
$$
   T(\delta_x(\lab  \eta_t ,\cdot \rab))(f) 
   	= \frac{1}{\sqrt{2\pi \lab \eta_t, \eta_t \rab}} 
		\exp \left( -\frac{1}{2} \frac{(i\lab \eta_t ,f \rab - x)^2}{\lab \eta_t,\eta_t \rab} -\frac{1}{2}\lab f,f\rab \right)
$$
for $f\in S(\R)$. Using 
$$
    \lab \eta_t,\eta_t\rab = \sigma^2 \int_0^t e^{2\lambda(s-t)}\, ds 
    		= \sigma^2 \frac{1-\exp(-2\lambda t)}{2\lambda} =: k
$$
the expectation is readily available by
$$
   T(\delta_x(\lab \eta_t, \cdot \rab)(0) = \frac{1}{\sqrt{2\pi k}} \exp \lb -\frac{x^2}{2k}\rb\;.
$$
\bigskip
\begin{proposition}[Expectation of occupation times]
Let $X_t$ be a centered Ornstein--Uhlenbeck process on the time interval $[0,T]$, where $T>0$. Let $[a,b]\subset\R$ be an interval, where $a=-\infty$ and $b=\infty$ are allowed. The expectation of the occupation time $M_{T,[a,b]}(X_t)$ is given by
\begin{eqnarray}
   \label{E:OccTime}
   \mathbb{E}_\mu\left( M_{T,[a,b]}(X_t) \right)& \nonumber\\
   =& \frac{1}{2} \int_0^T \erf\left( \frac{\alpha b}{1-\exp(-2\lambda t)}\right) -  \erf\left( \frac{\alpha a}{1-\exp(-2\lambda t)}\right) \, dt\;,
\end{eqnarray}
where $\alpha =\sqrt{\lambda}/\sigma$.
\end{proposition}
\bigskip
\begin{pf}
The occupation time, i.e.~the time the process spends in the space interval $[a,b]$ during the time $T$ is given by 
$$
    M_{T,[a,b]}(X_t) = \int_0^T \int_a^b \delta_0(X_t-x) \, dx \, dt\;.
$$
Interchanging integrations due to Theorem~\ref{T:Intcor}, we obtain the expectation of the occupation time of the Ornstein--Uhlenbeck process
\begin{eqnarray}
   \E_\mu(M_{T,[a,b]}(X_t)) &=& \int_0^T \int_a^b \E_{\mu}(\delta_0(X_t-x))\, dx\, dt \nonumber\\
   	&=&\int_0^T \frac{1}{\sqrt{2\pi k}} \int_a^b \exp \lb -\frac{x^2}{2k}\rb \, dx \, dt \nonumber \\
	&=& \frac{1}{2} \int_0^T \erf\lb \frac{b}{\sqrt{2k}}\rb-\erf\lb \frac{a}{\sqrt{2k}}\rb \, dt \nonumber
\end{eqnarray}
$$
	=\frac{1}{2} \int_0^T \erf\lb \frac{\alpha}{\sqrt{1-\exp(-2\lambda t)}}b\rb-\erf\lb \frac{\alpha}{\sqrt{1-\exp(-2\lambda t)}}a\rb \, dt,
$$
with $\alpha=\sqrt{\lambda}/\sigma$.
\end{pf}
\section{Numerics}
\label{S:Numerics}
\subsection{Estimation of the expected occupation time by Monte--Carlo methods}
The expected occupation time of the Ornstein--Uhlenbeck process $X_t$ defined via the SDE~(\ref{E:SDE_OU_SIM}) can be computed using Eqn.~(\ref{E:OccTime}). Alternatively, one can also compute the occupation time using a Monte--Carlo simulation of the underlying process. We generate $N$ of sample paths of the Ornstein--Uhlenbeck process with a fixed parameter set and compute the sample occupation time $\tilde{M}_{T,[a,b]}(\cdot)$ for each path.  As in the basic idea of the Monte--Carlo simulation, the sample average of the occupation time serves as an estimator for the expectation value. If large numbers of samples are considered, the estimator yields a better approximation. In the sequel, we shortly outline the numerical approximation of a stochastic process like the Ornstein--Uhlenbeck process~(\ref{E:SDE_OU_SIM}).

Consider a general non--autonomous stochastic differential equation
\begin{equation}
   \label{eqn:GENOUprocess}
   dX_t =f(t,X_t) dt + g(t,X_t) dW_t \;, \quad X_0=x_0\in \R
\end{equation}
defined in the time interval $[0,T]$, where $f,g:[0,T]\times\R\to \R$ and $W_t$ is a standard Wiener process. Under mild conditions the solution of (\ref{eqn:GENOUprocess}) has the following form
\begin{equation}
   \label{eq:intergralFormofSDE}
   X_t = X_0 + \int^t_0{f(s,X_s) ds} + \int^t_0 g(s,X_s) dW_s,\quad 0\leq t\leq T\;.
\end{equation}
Note that the solution $X_t$ is a random variable for each $t$. For details on the existence and uniqueness of solutions to (\ref{eqn:GENOUprocess}), we refer to Ref.~\cite{Oksa07}.

\begin{center}
\begin{figure}[h]
\includegraphics[width=\textwidth]{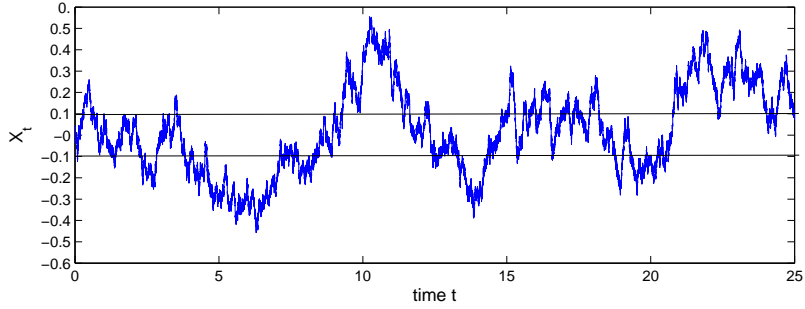} 
\caption{A sample path of the Ornstein--Uhlenbeck process~(\ref{E:SDE_OU_SIM}) with parameters $\lambda=0.5$ and $\sigma =0.25$ on the time interval $[0,25]$.} \label{fig:ouprocess}
\end{figure}
\end{center}

This solution can be numerically estimated by using the Euler--Maruyama method.
We discretize the interval $[0,T]$ using a time step $\Delta t = T/L$ for some positive integer $L$ and introduce discrete time points $\tau_j = j\Delta t$ for $j=1,2,...,L$. Let $\tilde{X}_j$ denote the numerical approximation of $X_{\tau_j}$. Further, we assume that the second integral on the right hand side of (\ref{eq:intergralFormofSDE}) is integrated using the It\^{o}--version of stochastic integrals. Then the Euler--Maruyama method reads as
\begin{equation}
   \label{eq:EMmethod}
   \tilde{X}_j = \tilde{X}_{j-1} +f(\tau_j,\tilde{X}_{j-1})\Delta t+g(\tau_j,\tilde{X}_{j-1})(W_{\tau_j} - W_{\tau_{j-1}}), \quad j = 1,2,...,L\;.
\end{equation}
Figure~\ref{fig:ouprocess} shows a sample path of the Ornstein--Uhlenbeck process computed using the Euler--Maruyama method. Using the discrete version $\tilde{X}_j$ of the process, we can easily calculate the sample occupation time $\tilde{M}_{T,[a,b]}(\tilde{X}_j)$.

\begin{remark}
The accuracy of the numerical solution to the SDE can be measured in two ways, namely strong and weak convergence. Strong convergence measures the accuracy on the basis of individual realizations. The weak convergence measures the accuracy of numerical methods to SDEs in case where the goal is to ascertain the probability distribution. For example, the Euler--Maruyama method has strong order of convergence $\gamma=\frac{1}{2}$. For more details, see Refs.~\cite{higham01,kloeden92}.
\end{remark}

\subsection{Direct computation of the expected occupation time}
To compute the expectation of the occupation times given by~(\ref{E:OccTime}), we have to evaluate integrals of the type
\begin{eqnarray}
   \label{E:erfint}
   \int_0^T \erf\lb \frac{C}{\sqrt{1-\exp(-2\lambda t)}}\rb\, dt 
   	&=& \frac{C^2}{\lambda} \int_{s_0}^\infty \frac{\erf(s)}{s(s^2-C^2)}\, ds \nonumber \\
	&:=& \frac{C^2}{\lambda} g(C,\lambda, T) 
\end{eqnarray}
where $s_0=C/\sqrt{1-\exp(-2\lambda T)}$. Note that $s_0>C$ and $s_0\to C$ for $\lambda T\to \infty$. Hence, in the limit $\lambda T\to \infty$, e.g.~for $\lambda=\mathcal{O}(1)$ and $T\to \infty$, the integral gets singular. This together with the unbounded domain of integration poses some numerical difficulties, which can be overcome by splitting the integral: the part close to the asymptotic singularity at $s=C$, and intermediate part and the part close to infinity. We introduce $s_1>s_0$ and $s_2>s_1$ and rewrite
$$
   g(C,\lambda, T) = \lb \int_{s_0}^{s_1} + \int_{s_1}^{s_2} + \int_{s_2}^{\infty} \rb 
   						\frac{\erf(s)}{s(s^2-C^2)}\, ds := I_1 + I_2 + I_3\;.
$$
The second integral $I_2$ does not pose any numerical difficulties and can easily be computed using Simpson's rule. However, we have to take care about the first and the third part.

 In $I_1$, we replace the error function by its quadratic Taylor polynomial 
$$
   T_2(s) = \erf(s_0) + \frac{2e^{-s_0^2}}{\sqrt{\pi}}\lsb (s-s_0)- s_0(s-s_0)^2\rsb
$$
at $s_0$, and get
\begin{eqnarray}
    I_1 &\sim \erf(s_0) \int_{s_0}^{s_1} \frac{ds}{s(s^2-C^2)} 
    	+ \frac{2e^{-s_0^2}}{\sqrt{\pi}} \int_{s_0}^{s_1} \frac{(s-s_0)-s_0(s-s_0^2)}{s(s^2-C^2)}\, ds \\
    	& = \frac{1}{2C^2} \lb \erf(s_0)-\frac{2(s_0^3+s_0)e^{-s_0^2}}{\sqrt{\pi}}\rb \nonumber
					\ln \frac{1-C^2/s_1^2}{1-C^2/s_0^2} \\
			& \qquad  + \frac{2s_0^2+1}{C\sqrt{\pi}} e^{-s_0^2} \ln \frac{(s_1-c)(s_0+C)}{(s_1+C)(s_0-C)}
		+ \frac{s_0^2}{\sqrt{\pi}} e^{-s_0^2} \ln\frac{s_1^2-C^2}{s_0^2-C^2}\;.
\end{eqnarray}
The choice of $s_1$ depends upon the desired accuracy of the above approximation.

\begin{lemma} \label{LM:lemmaIntEstim}
Choosing $s_1=s_0+\sqrt[4]{s_0(s_0^2-C^2)\varepsilon}$ for a given tolerance $\varepsilon>0$, we get an approximation error $\le \varepsilon$.
\end{lemma}
\bigskip
\begin{pf}
\begin{eqnarray}
   | \int_{s_0}^{s_1} \frac{\erf(s)-T_2(s)}{s(s^2-C^2)}\, ds | &\le& \max\limits_{s_0\le s\le s_1} |\erf(s)-T_2(s)| \cdot \int_{s_0}^{s_1} \frac{ds}{s(s^2-C^2)} \nonumber\\
																																&\le& \max\limits_{s\in \R} |\frac{d^3 \erf}{ds^3}(s)| \cdot \frac{d^3}{6} \cdot \frac{d}{s_0(s_0^2-C^2)} \nonumber,
\end{eqnarray}
where we introduce $d=s_1-s_0$ and estimate the integral using the mean value theorem. Hence
\begin{eqnarray}   
   | \int_{s_0}^{s_1} \frac{\erf(s)-T_2(s)}{s(s^2-C^2)}\, ds |
   	& \le \frac{2d^4}{3\sqrt{\pi}\, s_0(s_0^2-C^2)} \nonumber
\end{eqnarray}
For the given choice of $s_1$, i.e.~$d=\sqrt[4]{s_0(s_0^2-C^2)}$ we get
$$
   | \int_{s_0}^{s_1} \frac{\erf(s)-T_2(s)}{s(s^2-C^2)}\, ds | \le 0.376\varepsilon
$$
\end{pf}

\bigskip

In the third part of the integral, we replace the error--function by its limit $\erf(\infty)=1$ and get
$$
   \int_{s_2}^{\infty} \frac{\erf(s)}{s(s^2-C^2)}\, ds 
   					\approx \frac{1}{2C^2}\ln \frac{s_2^2}{s_2^2-C^2}\;.
$$
Choosing $s_2>10$ yields an approximation of the error function of less than $10^{-44}$.

\begin{center}
\begin{figure}
\includegraphics[width=\textwidth]{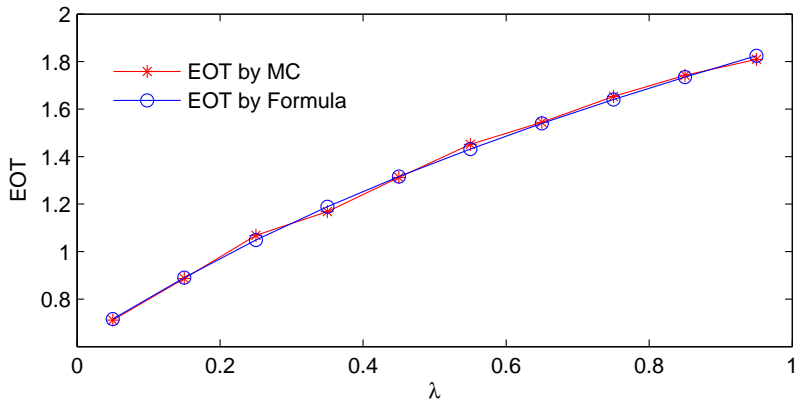} 
\caption{Expected occupation times on the interval $[a,b]=[-0.1, 0.1]$, time horizon $T=16$ and parameters $\sigma=1$ and $\lambda\in(0,1)$. Computations are carried out using either the Monte--Carlo method or Eqn.~(\ref{E:OccTime}).} 
\label{fig:optFordifLamSig_1}
\end{figure}
\end{center}

\noindent Figure \ref{fig:optFordifLamSig_1} shows the expec\-ted occupation times obtained using either Eqn.~(\ref{E:OccTime}) or Monte-Carlo simulations. In the setting that we have shown here, both computational methods yield indistinguishable results.

\begin{remark}
For large values of $2\lambda T$, the function $g(C,\lambda, T)$ in~(\ref{E:erfint}) gets hard to evaluate numerically. For $2\lambda T>37$, we obtain $e^{-2\lambda T}<10^{-16}$; the usual machine precision. Hence, we limit ourselves to $2\lambda T<37$. Furthermore, in case of $\lambda T \to \infty$, we obtain $s_0 \to C$ as well as $s_1 \to C $. Therefore the splitting of the integration introduced above does not resolve the problem with singularity at the bounds of the integration interval.
\end{remark}

\subsection{Parameter estimation}
\noindent Estimating the parameters of the Ornstein--Uhlenbeck process given by~(\ref{E:SDE_OU_SIM}) based on Eqn.~(\ref{E:OccTime}) for the occupation time is the main objective of this paper. Therefore, we formulate an optimization problem which we use to estimate the parameters.

Let $X_{\lambda,\sigma}(t)$ denote the Ornstein--Uhlenbeck process with parameters $\lambda$ and $\sigma$ as in Eqn.~(\ref{E:SDE_OU_SIM}). Let a spatial interval $[a,b]$ be given fixed. Then we define $G_{T,[a,b]} (\lambda, \sigma):= \mathbb{E}_\mu \left( M_{T,[a,b]}(X_{\lambda,\sigma}(t)) \right)$ as the expected occupation time of the Ornstein--Uhlenbeck process for time horizon $T$.

\medskip

\emph{Problem: Given different time horizons $T_i$ and intervals $[a_j,b_j]$, corresponding occupation times $G_{i,j}$ for $i=1,\dots,n$ and $j=1,\dots,m$, determine the parameters $(\lambda,\sigma)$ such that the deviation
\begin{equation}
   R(\lambda,\sigma):=\sum_{i=1}^n \sum_{j=1}^m\left(G_{T_i,[a_j,b_j]}(\lambda,\sigma)-G_{i,j}\right)^2
\end{equation}
is minimal.}

\bigskip

To solve this optimization problem, we apply a standard method from numerical analysis. Here, we have used the simplex search method implemented in \textsc{Matlab} as the function \textit{fminsearch}, see Ref.~\cite{matlabCAT07}.
As a stopping exit for the optimization, we use a tolerance of $10^{-5}$ between successive iterations.
To demonstrate the parameter estimation procedure, consider the following situation. Let $T_1 =10$ and $T_2=12$ be the time horizons $T_i$ for $i = 1,2$ and $[-0.25,0.25]$, 
$[-0.5,0.5]$, $[-0.75,0.75]$ $[-1.0,1.0]$ be the intervals $[a_j,b_j]$ for $j=1,2,3,4$, which use to calculate the corresponding data $G_{i,j}$ to the above mentioned optimization problem. Using either the direct equation~(\ref{E:OccTime}) or the Monte--Carlo method, we compute occupation times $G_{i,j}$ for the parameters $(\lambda,\sigma)=(0.15,0.90)$.  Now we initiate the minimization procedure providing the corresponding $G_{i,j}$s for the both situations separately. The resulting estimated values are $(\lambda^*,\sigma^*)=(0.150025,0.90004)$ in the case $G_{i,j}$ are computed via direct equation~(\ref{E:OccTime}) and $(\lambda^*,\sigma^*)=(0.1330915,0.8764664)$ in the case $G_{i,j}$ are computed by the Monte--Carlo methods. The following table lists more numerical findings that are estimated correspond to different setting.

\begin{center}
\begin{tabular}{cc|cc|cc}
\multicolumn{2}{c|}{true parameters} &
\multicolumn{2}{c|}{recovered from~(\ref{E:OccTime})} &
\multicolumn{2}{c}{recovered from MC}\\
$\lambda$ & $\sigma$ & $\lambda^\ast$ & $\sigma^\ast$ & $\lambda^\ast$ & $\sigma^\ast$ \\
\hline
0.25 & 0.75 & 0.250013& 0.750024 & 0.275525 & 0.779474\\
0.50 & 0.50 & 0.500012 & 0.499878 & 0.549382 & 0.519195\\
0.75 & 1.25 & 0.750022 & 1.249987 & 0.720177 & 1.225472\\
1.00 & 2.00 & 0.999917 & 2.000018 & 1.026446 & 2.026084\\
1.25 & 2.50 & 1.250011 & 2.500011 & 1.309042 & 2.554252\\
\end{tabular}
\end{center}
The parameters recoverd from occupation times generated using~(\ref{E:OccTime}) (columns 3 and 4) agree better than those recoverd from occupations times generated with the help of Monte--Carlo simulations (last two columns). This is not surprising, since we used the same underlying equation to generate and to recover the parameters. However, also for the parameters covered from the Monte--Carlo simulations we have a difference of about 10$\%$ between the true and the recovered data. This is for most applications a sufficient accuracy. Nevertheless, increasing  the number of samples in Monte--Carlo simulation we can improve the accuracy of the recomputed parameters.

\section{Conclusion}
\label{S:Conclusion}
We derived an equation to compute directly the expected occupation time of the centered Ornstein--Uhlenbeck process. This allows to identify the parameters of the Ornstein--Uhlenbeck process for available occupation times via a standard least squares minimization. To test our method, we generated occupation times and recovered the parameters with the above mentioned procedure. Within the range of our numerical experiments, we found very good agreement. This gives hope to be able to estimate parameters in industrial fleece production processes from measurable quantities like the mass per area. However, to get closer to the industrial applications we have to extend our method to the $2d$--case and more involved processes than the standard Ornstein--Uhlenbeck process.

\section*{Acknowledgements}
\noindent W.~Bock would like to thank the Department of Mathematics for the financial support. The DAAD (German Academic Exchange Service) is gratefully acknowledged for providing a scholarship to U.~P.~Liyanage.

\end{document}